\documentstyle{article}
\begin{document}
\title{{Non-classical Solution to Hessian   
Equation   from Cartan Isoparametric Cubic}} \author{{Nikolai Nadirashvili\thanks{CMI, 39, rue F. Joliot-Curie, 13453
Marseille  FRANCE, nicolas@cmi.univ-mrs.fr},\hskip .4 cm Vladimir Tkachev\thanks
{Department of mathematics, KTH
SE-100 44 Stockholm SWEDEN, tkatchev@kth.se},\hskip .4 cm Serge
Vl\u adu\c t\thanks{IML, Luminy, case 907, 13288 Marseille Cedex
FRANCE, vladut@iml.univ-mrs.fr} }}

\date{}
\maketitle

\def\n{\hfill\break} \def\al{\alpha} \def\be{\beta} \def\ga{\gamma} \def\Ga{\Gamma}
\def\om{\omega} \def\Om{\Omega} \def\ka{\kappa} \def\lm{\lambda} \def\Lm{\Lambda}
\def\dl{\delta} \def\Dl{\Delta} \def\vph{\varphi} \def\vep{\varepsilon} \def\th{\theta}
\def\Th{\Theta} \def\vth{\vartheta} \def\sg{\sigma} \def\Sg{\Sigma}
\def\bendproof{$\hfill \blacksquare$} \def\wendproof{$\hfill \square$}
\def\holim{\mathop{\rm holim}} \def\span{{\rm span}} \def\mod{{\rm mod}}
\def\rank{{\rm rank}} \def\bsl{{\backslash}}
\def\il{\int\limits} \def\pt{{\partial}} \def\lra{{\longrightarrow}}

\section{Introduction}
\bigskip

This paper shows that one can use certain specific minimal cubic cone, namely, 
the Cartan isoparametric eigencubic [C] to construct non-smooth solution to Hessian 
fully nonlinear second-order elliptic equation.

More precisely  we study a class of fully nonlinear second-order elliptic equations
 of the form
$$F(D^2u)=0\leqno(1.1)$$
defined in a domain of ${\bf R}^n$. Here $D^2u$ denotes the
Hessian of the function $u$. We assume that
$F$ is a Lipschitz  function defined on the space $ S^2({\bf R}^n)$ 
of ${n\times n}$ symmetric matrices  satisfying the uniform ellipticity condition,
 i.e. there exists a constant $C=C(F)\ge 1$ (called an {\it ellipticity
constant\/}) such that 
$$C^{-1}||N||\le F(M+N)-F(M) \le C||N||\;
\leqno(1.2)$$ 
for any non-negative definite symmetric matrix $N$; if
$F\in C^1(S^2({\bf R}^n))$ then this condition is equivalent to
$${1\over C'}|\xi|^2\le F_{u_{ij}}\xi_i\xi_j\le C' |\xi |^2\;,
\forall\xi\in {\bf R}^n\;.\leqno(1.2')$$
 Here, $u_{ij}$ denotes the partial derivative
$\pt^2 u/\pt x_i\pt x_j$. A function $u$ is called a {\it
classical\/} solution of (1) if $u\in C^2(\Om)$ and $u$ satisfies
(1.1).  Actually, any classical solution of (1.1) is a smooth
($C^{\alpha +3}$) solution, provided that $F$ is a smooth
$(C^\alpha )$ function of its arguments.

For a matrix $S \in   S^2({\bf R}^n)$  we denote by $\lambda(S)=\{
\lambda_i : \lambda_1\leq...\leq\lambda_n\}
 \in {\bf R}^n$  the (ordered) set  of eigenvalues of the matrix $S$.
Equation (1.1) is called a {\em Hessian equation} ([T1],[T2] cf. [CNS])
 if the function $F(S)$ depends only on the eigenvalues $\lambda(S)$ of the matrix $S$, i.e., if
 $$F(S)=f(\lambda(S)),$$
 for some function $f$  on ${\bf R}^n$ invariant under  permutations of
 the coordinates.

 In other words the equation (1.1) is called Hessian if it is invariant under
 the action of the group
 $O(n)$ on $S^2({\bf R}^n)$:
 $$\forall O\in O(n),\; F({^t O}\cdot S\cdot O)=F(S) \;.\leqno(1.3) $$
 The Hessian invariance relation (1.3) implies the following:

\medskip
 (a) $F$ is a smooth (real-analytic) function of its arguments if and only if $f$ is
a smooth (real-analytic) function.

\medskip
 (b) Inequalities (1.2) are equivalent to the inequalities
 $${\mu\over C_0} \leq { f ( \lambda_i+\mu)-f ( \lambda_i) } \leq C_0 \mu,
 \; \forall  \mu\ge 0,$$
 $\forall  i=1,...,n$, for some positive constant $C_0$.

\medskip
 (c) $F$ is a concave function if and only if $f$ is concave.

\medskip
 Well known examples of the Hessian equations are Laplace, Monge-Amp\`ere,
Bellman, Isaacs and Special Lagrangian equations.

 \medskip
 Bellman and Isaacs equations appear in the theory of controlled diffusion processes, see [F]. 
 Both are fully nonlinear uniformly elliptic equations of the form (1.1). The Bellman equation
 is concave in $D^2u \in  S^2({\bf R}^n)$ variables. However, Isaacs operators are, in  general,
 neither concave nor convex. In a simple homogeneous form the Isaacs equation can be written
 as follows:
$$F(D^2u)=\sup_b \inf_a L_{ab}u =0, \leqno (1.4) $$
 where $L_{ab}$ is a family of linear uniformly elliptic operators of type 
 $$L=  \sum a_{ij} {\partial^2 \over  \partial x_i \partial x_j } \leqno (1.5)$$
  with an ellipticity  constant $C>0$ which depends on two parameters $a,b$.

\medskip
Consider the Dirichlet problem
$$\cases{F(D^2u)=0 &in $\Om$\cr
u=\vph &on $\pt\Om\;,$\cr}\leqno(1.6)$$ where  $\Omega \subset {\bf R}^n$ is a
 bounded domain with smooth boundary $\partial \Omega$
and $\vph$ is a continuous function on $\pt\Om$.

We are interested in the problem of existence and regularity of
solutions to the Dirichlet problem (1.6) for Hessian equations and Isaacs equation.
The problem (1.6) has always a unique viscosity (weak)
solution for  fully nonlinear elliptic equations (not necessarily
Hessian equations). The viscosity solutions  satisfy the equation
(1.1) in a weak sense, and the best known interior regularity
([C],[CC],[T3]) for them is $C^{1,\epsilon }$ for some $\epsilon
> 0$. For more details see [CC], [CIL]. Note, however, that viscosity solutions 
are $C^{2,\epsilon }$-regular  almost everywhere; in fact, it is true on the complement of  a close set of strictly less Hausdorff dimension [ASS].
Until recently it remained
unclear whether non-smooth viscosity solutions exist. In the recent papers [NV1], [NV2],  [NV3], [NV4] two of the present authors first proved   the existence   of non-classical viscosity solutions to a fully nonlinear elliptic equation, and then of singular solutions to Hessian  uniformly elliptic equation in all dimensions beginning from 12. Those papers   use the functions 
$$w_{12,\delta}(x)={P_{12}(x)\over |x|^{\delta }},\;w_{24,\delta}(x)= {P_{24} (x)\over  |x|^{\delta }},\:
\delta\in [1,2[,
$$ 
with $P_{12}(x),P_{24}(x)$ being cubic forms as follows:
 $$P_{12}(x)=Re (q_1q_2q_3),\; x=(q_1,q_2,q_3)\in {\bf H}^3={\bf R}^{12},$$
  $  {\bf H}$ being Hamiltonian quaternions,
$$ P_{24}(x)={Re((o_1\cdot o_2)\cdot o_3)}={Re(o_1\cdot(o_2\cdot o_3))},\; x=(o_1,o_2,o_3)\in {\mathcal{O}}^3={\bf
R}^{24}
$$   $\mathcal{O}$ being the  algebra of Caley octonions.
\smallskip

 As it was noted by the second author (V.T.), these are (most symmetric) examples of so-called {\em radial eigencubics}, which  define minimal cubic cones. Since the family of such cubics  contains especially interesting isoparametric  Cartan cubics  [C]  in dimensions 5, 8, 14 and 26 admitting large automorphism groups, the last two  being intimately connected with $P_{12}(x)$ and $P_{24}(x)$, it is but natural to try these Cartan cubics as numerators of tentative non-classical solutions to Hessian  uniformly elliptic equation in their respective dimensions. Our  main goal in this paper is to  show that   in 5 dimensions  it really works at least for $C^{1,1}$-solutions, and to prove

\medskip

{\bf Theorem 1.1.}
\medskip
{\it The function
$$w_5(x)={P_{5} (x)\over |x|}
$$ is a viscosity solution  to  a uniformly elliptic Hessian equation
  $(1.1)$ in a unit ball $B\subset {\bf
R}^{5}$  for  the isoparametric Cartan cubic form in
$x=(x_1,x_2,z_1,z_2,z_3)\in {\bf R}^{5}$
$$ P_{5}(x)=x_1^3+\frac{3 x_1}2\left(z_1^2 + z_2^2-2 z_3^2-2x_2^2\right)+\frac{3\sqrt 3}2\left(x_2z_1^2-x_2z_2^2 + 2z_1z_2z_3\right).$$
 }

\medskip 
\noindent At the time  of writing it is not at all clear that the same is true for the function 
$w_{5,\delta}(x)=P_{5} (x)/ |x|^{\delta }$ for $\delta>1$, see Remark 4.1 below, and thus 
the question on the optimality of  the interior $C^{1,\alpha}$-regularity of viscosity solutions to fully nonlinear 
equations is open   in dimensions up to 12.
\medskip

However, the method of [NV2] permits to  construct singular solutions in ten dimensions:

\medskip
{\bf Corollary 1.1.}
\medskip

{\it There exist $\delta>0, M>0$ such that  the  homogeneous order $(2-2\delta)$ function 
$$u_{10,\delta,M}(x,y)={w_5(x)+w_5(y)+M(|x|^2-|y|^2)\over (|x|^2+|y|^2)^{\delta } }$$ 
in the unit ball $B\subset {\bf R}^{10}$   is a viscosity solution to a uniformly elliptic equation $(1.1)$.}
\medskip

For a proof it is sufficient just to repeat the argument of [NV2] which gives the result for
 $\delta= 10^{-6},  M=100$.

\medskip
As in [NV3] we get  also that $w_5$  is a viscosity solution  to
a uniformly elliptic Isaacs equation:
\bigskip

{\bf Corollary 1.2.}

\medskip
{\it The function
$$w_5(x)=P_{5} (x)/ |x| 
$$ is a viscosity solution  to  a uniformly elliptic  Isaacs equation
  $(1.4)$ in a unit ball $B\subset {\bf
R}^{5}$.}
\medskip

{\em Remark 1.1.} One could hope that using a minimal cubic cone in 
4 dimensions, namely, the   Lawson cubic [L] (essentally unique by [P])
 $$P_4(x)=x_3(x_1^2-x_2^2)+2x_1x_2x_4,$$ 
it is possible to costruct a non-smooth solution to a Hessian uniformly elliptic equation in four dimensions, 
but it does not work. Namely, the function $$w_4(x)=P_4 (x)/ |x| $$
{\em is not} a solution to such  equation, since it does not verify the
 conditions of Lemma 2.1 (the corresponding matrix family is not hyperbolic).

\medskip
 The rest of the paper is organized as follows: in Section 2 we recall some necessary preliminary results from 
[NV3], then we recall some facts about radial eigencubics and especially the Cartan cubic in Section 3 and we 
prove our main results in Section 4.\smallskip
  
One should note that similar results are valid for three other Cartan cubics and that all of them can be used in some other similar applications. These results which depend on the theory of Jordan algebras will be exposed elsewhere.

\section{Preliminary results }

\medskip
   Let $w=w_n$ be an odd homogeneous  function
of order 2, defined on a unit ball
$B =B_1\subset {\bf R}^n$ and smooth in $B \setminus\{0\}$. Then
the Hessian of $w$ is homogeneous of order 0. 

We want to give a criterion for $w$ to be a solution of a uniformly elliptic Hessian  equation
or  a uniformly elliptic Isaacs equation. To do this, recall that  a family 
${\mathcal A}\subset S^2({\bf R}^n)$ of
symmetric matrices $A$ is  called {\em uniformly hyperbolic} if there exists  a  constant $M>1$ such that 
$$\frac{1}{ M} < -\frac{\lambda_1(A)}{\lambda_n (A)} < M$$
for any $A\in{\mathcal A},\; \lambda_1(A)\geq...\geq\lambda_n(A)$ being the eigenvalues of $A$.\medskip

One can reformulate some results from [NV3] (namely, Lemma 2.1,  six final lines of Section 4, 
Lemmas 5.1 and 5.2) as follows, in our special case $\delta=1$:

\medskip
{\bf Lemma 2.1.} 

\medskip
{\em Set for $x,y\in S^{n-1}$ and for an orthogonal matrix $O\in  O(n)$,
 $$M(x,y,O):=D^2w(x)-{^t O}D^2w(y) O.$$
Suppose that the family $${\mathcal M}:=\{M(x,y,O):M(x,y,O)\neq 0, x\neq y, x\neq 0,y\neq 0,O\in  O(n)\}\subset S^2({\bf R}^n)$$ is uniformly hyperbolic. Then $w$ is a solution to a uniformly elliptic Hessian  equation
as well as to  a uniformly elliptic Isaacs equation.}\bigskip

We need also the following property of the eigenvalues $
\lambda_1 \ge\ldots\ge\lambda_{n}$ of real symmetric matrices 
 of order $n$ which is a classical result by Hermann  Weyl [W]:

\bigskip\noindent {\bf  Lemma 2.2.} \medskip

  {\em Let $ A\neq B$  be two real symmetric  matrices
with the eigenvalues $
\lambda_1\ge\lambda_2\ge\ldots\ge\lambda_{n} $ and $
\lambda'_1\ge\lambda'_2\ge\ldots\ge\lambda'_{n} $ respectively.
Then for the eigenvalues $
\Lambda_1\ge\Lambda_2\ge\ldots\ge\Lambda_{n} $ of the matrix $A-B$
we have}
$$ \Lambda_1\ge\max_{i=1,\cdots, n}(\lambda_i-\lambda'_i),
 \;\;\Lambda_n\le\min_{i=1,\cdots, n}(\lambda_i-\lambda'_i). $$

\medskip

\section{Radial eigencubics }
Let us recall the Cartan  cubic form $P_5(x)$ which is 
closely  related with  real algebraic minimal cones, that is, homogenous polynomial solutions $u$ to the minimal surface equation
$$(1 + |\nabla u|^2)\Delta u -\sum u_{ij}u_{i}u_{j}=0 .$$

According to  Hsiang [H], the study of those is equivalent to  classifying  polynomial
solutions $f = f(x_1, . . . , x_n)\in {\bf R}[x_1, . . . , x_n]$ of the following congruence:
$$L(f) =0 \:(\mod f), \leqno (3.1)$$
 $$L(f): = |\nabla f|^2\Delta f -\sum f_{ij}f_{i}f_{j}$$
being the normalized mean curvature operator. This condition means that
the zero-locus $f^{-1}(0)$ has zero mean curvature everywhere where the gradient $\nabla f\neq 0$.
A non-zero polynomial solutions of this congruence is called an {\em eigenfunction} of $L.$
 The ratio $L(f)/f$ (a polynomial in $x$) is called the weight of an eigenfunction $f.$ An eigenfunction $f$
which is a cubic homogenous   is  called an  {\em  eigencubic}.

 Among them the most interesting  are solutions of the following
non-linear equation:
$$L(f) = \lambda |x|^2f,\; \lambda\in {\bf R}, \leqno (3.2)$$
which are called {\em radial eigencubics}.  In [H], Hisang posed the problem to
 determine all solutions of (3.2) up to a congruence in ${\bf R}^n$ (for any degree).

This classification for  radial eigencubics is almost completed in [Tk1,Tk2], namely,
 any radial eigencubic is either a member of the infinite family of eigencubics of
 Clifford type completely classified in [Tk1], or belongs to one
of exceptional families, the number of these  lying between 13 and 19.

The cubic forms $P_{12}(x)$ and $P_{24}(x)$ belong to the Clifford family; the Cartan polynomial
$P_{5}(x)$ is the first, that is of least dimension,  in the list of  exceptional radial eigencubics. Moreover,
it is an isoparametric polynomial, that is satisfies the M\"unzner  system [M]:
$$|\nabla f|^2=9|x|^4, \;\Delta f=0, $$
expressing the fact that all principal cuvatures of $f^{-1}(0)\bigcap S^4$ are constant 
(and different). Since $L(f) =-54 |x|^2 f$, $P_{5}(x)$ is  a  radial eigencubic as well.

The form  $P_{5}(x)$ admits a three-dimensional automorphism group. Indeed, one easily 
verifies that the orthogonal trasformations 

$$ A_1(\phi):= { 1\over 2}\left(%
\begin{array}{ccccc}
   {3\cos(\phi)^2-1 }& { \sqrt 3\sin(\phi)^2 }&{ \sqrt 3\sin(2\phi) } &0&0\\
{ \sqrt 3\sin(\phi)^2 }&{1+\cos(\phi)^2 }&{-\sin(2\phi) } &0&\;0\\
 {- \sqrt 3\sin(2\phi) }&{\sin(2\phi) }&2\cos(2\phi) &0&0\\
  0& 0&0&2\cos(\phi)&2\sin(\phi) \\
   0&0&0&-2\sin(\phi)&2\cos(\phi)\\
 
\end{array}%
\right) $$

$$ A_2(\psi):= \left(%
\begin{array}{ccccc}
   1& 0& 0& 0& 0\\
0&\cos(2\psi)& 0& -\sin(2\psi)& 0\\
0&0& \cos(\psi)&0&-\sin(\psi)\\
0&\sin(2\psi)& 0&\cos(2\psi)& 0\\
0&0&\sin(\psi)& 0& \cos(\psi)\\
\end{array}%
\right) $$

$$ A_3(\theta):=  { 1\over 2}\left(%
\begin{array}{ccccc}
{3\cos(\theta)^2-1 }& {-\sqrt 3\sin(\theta)^2 }&0& 0&  {-\sqrt 3\sin(2\theta) } \\ 
{-\sqrt 3\sin(\theta)^2 }& {1+\cos(\theta)^2 }& 0& 0& {-\sin(2\theta) }\\
0&0&2\cos(\theta)&-2\sin(\theta)&0\\
0& 0& 2\sin(\theta)&2\cos(\theta)& 0\\
{\sqrt 3\sin(2\theta) }&  {\sin(2\theta) }& 0& 0&2\cos(2\theta)\\
\end{array}%
\right) $$

do not change the value of $P_5(x)$.

\medskip
Moreover, one easily gets

\medskip

{\bf Lemma 3.1.}
\medskip

 {\em Let $G_P$ be subgroup of $SO(5)$ generated by

 $$\{A_1(\phi),A_2(\psi), A_3(\theta):\:(\phi,\psi,\theta)\in {\bf R}^3\}.$$
 Then the orbit $G_PS^1$ of
 the circle
 $$S^1=\{(\cos(\chi),0,\sin(\chi),0,0):\chi\in {\bf R}\}\subset S^4  $$
under the natural action of  $G_P$ is the whole $S^4 . $ }

\medskip
{\em Proof.} Indeed, calculating  the differential of the action 
$$ (S^1)^4\longrightarrow S^4,\: (\phi,\psi,\theta,\chi)\mapsto
 (\cos(\chi),0,\sin(\chi),0,0){A_1(-\phi)}{A_2(-\psi)}{A_3(-\theta)}$$
at $(\phi,\psi,\theta,\chi)=(0,0,0,0)$ one sees that 
its rank is 4 which implies the sujectivity.

\section{Proofs}

Let $w_5=P_5(x)/|x|$. By Lemma 2.1  it is sufficient to prove the uniform hyperbolicity of the family
 $$M_5(x,y,O):=D^2w_5(x)-{^t O}D^2w_5(y) O.$$
\medskip

{\bf Proposition 4.1.}\medskip

 {\em Let $O\in O(5), x, y\in S^4$, $M_5(x,y,O)\neq 0$ and let 
$\Lambda_1\ge\ldots\ge\Lambda_5 $  be its eigenvalues.
Then $$\frac{1}{20}\le -\frac{\Lambda_1}{\Lambda_5}\le 20.  $$  }

{\em Proof.}  We beging with calculating the eigenvalues of $D^2w_5(x)$.

More precisely, we need \medskip

{\bf Lemma 4.1.} \medskip

{\em Let  $x\in S^4$, let $\lambda_1\ge\lambda_2\ge\ldots\ge\lambda_5 $  be the eigenvalues
of $D^2w_5(x)$, and let   $x\in G_P(p,0,q,0,0)$ with $p^2+q^2=1$.
Then $$ \lambda_1=\frac{p^3-6p+3 \sqrt{3(4-p^2)}}{2},\;\lambda_3= \frac{p^3+3p}{2},\; \lambda_5=\frac{p^3-6p-3 \sqrt{3(4-p^2)}}{2}. $$}

{\em Proof of Lemma 4.1.} Since $w_5$ is invariant under $G_P$, we can suppose that $x=(p,0,q,0,0)$.
Then $w_5(x)=\frac{p(3-p^2)}{2}$ and we get by a brute force calculation:
$$ D^2w_5(x):= \frac{ 1}{2}\left(%
\begin{array}{cc} M_1&0\\
0&M_2
  \\
\end{array}%
\right) $$ being a block matrix with
$$ M_1:= \left(%
\begin{array}{ccc}
   { p(1+2 p^2-3p^4)}&  {3\sqrt 3  p   ( p^2-1)}& {3 q(1-p^4)} \\ 
{3\sqrt 3  p   ( p^2-1)}& { p^3-15p}&  {3\sqrt 3  q   ( p^2+1)}\\
 {3 q(1-p^4)}& {3\sqrt 3  q   ( p^2+1)}& { p^3+3p^5} \\
  
\end{array}%
\right), \;$$
 $$
 M_2:= \left(%
\begin{array}{cc}
{ p^3+3p}& 6\sqrt 3 q \\
 6 \sqrt 3 q &  { p^3-15p}\\ 
\end{array}%
\right) $$
which gives for the characteristic polynomial $F(S) =F_1(S)\cdot F_2(S)$ where
$$ F_1(S)=\left(S-{3 p\over 2}-{ p^3\over 2}\right)\left(S^2+6 p S-p^3 S+{63 p^2\over 4}-3 p^4+
 {p^6 \over 4}-27\right)\:,$$ 
$$F_2(S)=\left(S^2+{15 p S\over 2}- {5p^3S\over 2} -{45 p^2\over 4}+ 
{15 p^4\over 2}-{5p^6\over 4}-9\right)$$
have  the roots

$$\lambda_1= \frac{p^3-6p+3 \sqrt{3(4-p^2)}}{2},\;\lambda_3= \frac{p^3+3p}{2},\; \lambda_5=\frac{p^3-6p-3 \sqrt{3(4-p^2)}}{2}\:, $$ 
$$\lambda_2= {5p^3-15p+3r\over 4},\; \lambda_4= {5p^3-15p-3r\over 4}\:$$
with $r :=\sqrt{5p^6-30p^4+45p^2+16}$.

One needs  only  to verify that indeed 
$$\lambda_1\ge\lambda_2\ge\lambda_3\ge\lambda_4\ge\lambda_5 $$
which is elementary. For example, let us verify for $p\in[-1,1]$ the inequality
$$\lambda_1=\frac{p^3-6p+3 \sqrt{3(4-p^2)}}{2}\ge\lambda_2= {5p^3-15p+3r\over 4},\; $$
which by symmetry gives $\lambda_4\ge\lambda_5 $ (the two remaining
inequalities being simpler).\smallskip

Indeed,

\noindent $ \lambda_1-\lambda_2 =3(r_1-r)/4$ with
$ r_1:=p-p^3+2\sqrt{3(4-p^2)}>0,$
$$r_1^2-r^2=4(1-p^2)(s_1+s_2) \ge 0$$
for $s_1:=p^4-6p^2+8=(4-p^2)(2-p^2)>0,\;s_2:=p\sqrt{3(4-p^2)}$
since
$$ s_1^2-s_2^2= (1-p)(4-p^2)(p^5+p^4-7p^3-7p^2+13p+16)\ge 0$$
because
$$p^5+p^4-7p^3-7p^2+13p+16=(p+1)(p^4-7p^2+13)+3\ge 3.$$

{\em End of proof.} Let now $y\in G_P(\bar p,0,\bar q,0,0)$. If $p=\bar p$ but $M_5(x,y,O)\neq 0$,
  the  trace $Tr(M_5(x,y,O))=8w_5(y)-8w_5(x)=0$ and the conclusion follows as for any traceless matrix in dimension 5. Let then $p>\bar p$;  by Lemma 2.2. one gets
 $$\Lambda_1\ge \lambda_2(p)- \lambda_2(\bar p)={(p-\bar p)(p^2+p\bar p+\bar p^2+3)\over 2}\ge
{3(p-\bar p)\over 2}, $$
 $$-\Lambda_5\ge \max\{\lambda_3(\bar p)-\lambda_3(p),\lambda_1(\bar p)-\lambda_1(p)\}\ge$$ $$\ge (p-\bar p)
\inf_{p\in[-1,1]}\max\{|\lambda'_1(p)|,|\lambda'_3(p)|\}=3(p-\bar p),$$
 $$0\le -Tr(M_5(x,y,O))=8w_5(x)-8w_5(y)=8(p-\bar p)(3-p^2-p\bar p-\bar p^2)\le 24 (p-\bar p).$$
Therefore,
$$-\Lambda_5\ge 4 \Lambda_1 $$
since $Tr(M_5(x,y,O))\le 0$ and
$$ 4 \Lambda_1+ \Lambda_5\ge Tr(M_5(x,y,O))\ge 24 (\bar  p-p), $$
 $$ - \Lambda_5\le 4 \Lambda_1 +24 (p-\bar p)\le 4 \Lambda_1+16 \Lambda_1=20\Lambda_1, $$
the case $p<\bar p$ being completely parallel which finishes the proof of our results.

\bigskip

{\em Remark 4.1.} One can   calculate the eigenvalues of $D^2(P_5(x)/|x|^{\delta})$ for $\delta>1$ as well
but then   the analogue of Lemma 4.1 does not hold, which makes impossible 
 applying the technique  of [NV3, Section 4]  to prove the hyperbolicity of the corresponding matrix family.
Thus the question whether  $P_5(x)/|x|^{\delta}$ is  for $\delta>1$ a solution of a uniformly
 elliptic (Hessian or not) equation remains open.

\bigskip
 \centerline{REFERENCES}

\medskip
 \noindent [ASS] S. N. Armstrong, L. Silvestre, C. K. Smart, {\em Partial regularity of solutions of fully nonlinear uniformly elliptic equations}, arXiv:1103.3677.

\medskip
 \noindent [C] L. Caffarelli,  {\it Interior a priory estimates for solutions
 of fully nonlinear equations}, Ann. Math. 130 (1989), 189--213.

\medskip
 \noindent [CC] L. Caffarelli, X. Cabre, {\it Fully Nonlinear Elliptic
Equations}, Amer. Math. Soc., Providence, R.I., 1995.

\medskip
 \noindent [CIL]  M.G. Crandall, H. Ishii, P-L. Lions, {\it User's
guide to viscosity solutions of second order partial differential
equations,} Bull. Amer. Math. Soc. (N.S.), 27(1) (1992), 1--67.

\medskip
 \noindent [CNS] L. Caffarelli, L. Nirenberg, J. Spruck, {\it The Dirichlet
 problem for nonlinear second order elliptic equations III. Functions
  of the eigenvalues of the Hessian, } Acta Math.
   155 (1985),  261--301.

  \medskip 
 \noindent [C]  \'E. Cartan, {\it Sur des familles remarquables d'hypersurfaces isoparam\'etriques dans les espaces sph\'eriques,} Math. Z. 45 (1939), 335-367.   

   \medskip
 \noindent [F] A. Friedman, {\it Differential games, } Pure and Applied Mathematics,
 vol. 25, John Wiley and Sons, New York, 1971.

\medskip
 \noindent [H] W. Y. Hsiang, {\it Remarks on closed minimal submanifolds in the standard Riemannian m-sphere,}
 J. Diff. Geom. 1(1967), 257-267.

\medskip
 \noindent [L]  H. B. Lawson, {\it Complete minimal surfaces in $S^3$}. Ann. of Math. 92 (1970), 335-374.

\medskip
 \noindent [M] H.-F. M\"unzner, {\it Isoparametrische Hyperfl\"achen in Sph\"aren.} Math. Ann. 251 (1980), 57-71.

\medskip
\noindent [NV1] N. Nadirashvili, S. Vl\u adu\c t, {\it
Nonclassical solutions of fully nonlinear elliptic equations,}
Geom. Func. An. 17 (2007), 1283--1296.

\medskip
\noindent [NV2] N. Nadirashvili, S. Vl\u adu\c t, {\it Singular Viscosity Solutions to Fully Nonlinear 
Elliptic Equations},  J. Math. Pures Appl., 89 (2008), 107-113.

\medskip
\noindent [NV3] N. Nadirashvili, S. Vl\u adu\c t, {\it Octonions and Singular Solutions of Hessian  Elliptic
Equations}, Geom. Func. An. 21 (2011), 483-498.

\medskip
\noindent [NV4] N. Nadirashvili, S. Vl\u adu\c t, {\it Singular solutions of Hessian fully nonlinear elliptic
equations}, Adv. Math., 228 (2011), 1718-1741.

\medskip 
 \noindent
[P] O.Perdomo, {\it Characterization of order 3 algebraic immersed minimal surfaces of $S^3$.} Geom. Dedicata, 129 (2007), 23-34.

\medskip 
 \noindent [Tk1] V.G. Tkachev, {\it On a classification of minimal cubic cones in ${\bf R}^n$}, arXiv:1009.5409.

\medskip 
 \noindent
[Tk2]   V.G. Tkachev, {\it Minimal cubic cones via Clifford algebras,} Compl. An. Oper. Th., 4 (2010), 
685-700. Preprint arXiv:1003.0215.

\medskip
\noindent [T1] N. Trudinger, {\it Weak solutions of Hessian
equations,} Comm. Part. Diff. Eq. 22 (1997),  1251--1261.

\medskip
\noindent [T2] N. Trudinger, {\it On the Dirichlet problem for
Hessian equations,} Acta Math. 175 (1995),  151--164.

\medskip
\noindent [T3] N. Trudinger, {\it H\"older gradient estimates for
fully nonlinear elliptic equations,} Proc. Roy. Soc. Edinburgh
Sect. A 108 (1988), 57--65.

\medskip
\noindent [We] G. Weyl, {\it Das asymptotische Verteilungsgezets
des Eigenwerte lineare partieller Differentialgleichungen,} Math.
Ann. 71 (1912),   441--479.
\end{document}